\documentclass{amsart}
\usepackage[dvips]{graphicx}
\usepackage{amscd}
\usepackage{amsmath}
\usepackage{amsxtra}
\usepackage{amsfonts}
\usepackage{amssymb}
\newtheorem{theorem}{Theorem}[section]

\newtheorem{lemma}[theorem]{Lemma}

\theoremstyle{definition}

\theoremstyle{remark}

\renewcommand{\theclaim}{\textup{\theclaim}}

\numberwithin{equation}{section}

\def\openone

{\mathchoice

{\hbox{\upshape \small1\kern-3.3pt\normalsize1}}

{\hbox{\upshape \small1\kern-3.3pt\normalsize1}}

{\hbox{\upshape \tiny1\kern-2.3pt\SMALL1}}

{\hbox{\upshape \Tiny1\kern-2pt\tiny1}}}

\makeatletter

\newbox\ipbox

\newcommand{\diracb}[1]{\left\langle #1\mathrel{\mathchoice

{\setbox\ipbox=\hbox{$\displaystyle \left\langle\mathstrut #1\right.$}

\vrule height\ht\ipbox width0.25pt depth\dp\ipbox}

{\setbox\ipbox=\hbox{$\textstyle \left\langle\mathstrut #1\right.$}

\vrule height\ht\ipbox width0.25pt depth\dp\ipbox}

{\setbox\ipbox=\hbox{$\scriptstyle \left\langle\mathstrut #1\right.$}

\vrule height\ht\ipbox width0.25pt depth\dp\ipbox}

{\setbox\ipbox=\hbox{$\scriptscriptstyle \left\langle\mathstrut #1\right.$}

\vrule height\ht\ipbox width0.25pt depth\dp\ipbox}

}\right. }

\newcommand{\dirack}[1]{\left. \mathrel{\mathchoice

{\setbox\ipbox=\hbox{$\displaystyle \left.\mathstrut #1\right\rangle$}

\vrule height\ht\ipbox width0.25pt depth\dp\ipbox}

{\setbox\ipbox=\hbox{$\textstyle \left.\mathstrut #1\right\rangle$}

\vrule height\ht\ipbox width0.25pt depth\dp\ipbox}

{\setbox\ipbox=\hbox{$\scriptstyle \left.\mathstrut #1\right\rangle$}

\vrule height\ht\ipbox width0.25pt depth\dp\ipbox}

{\setbox\ipbox=\hbox{$\scriptscriptstyle \left.\mathstrut #1\right\rangle$}

\vrule height\ht\ipbox width0.25pt depth\dp\ipbox}

} #1\right\rangle}

\newcommand{\linft}{L^{\infty}\left(  \mathbb{T}\right)}

\newcommand{\loner}{L^{1}\left(\mathbb{R}\right)}

\newcommand{\ltwor}{L^{2}\left(\mathbb{R}\right)}

\newcommand{\linfr}{L^{\infty}\left(\mathbb{R}\right)}

\input cyracc.def

\begin{document}
\title[The wavelet Galerkin operator]{The Spectrum of the Wavelet Galerkin Operator}
\author{Dorin Ervin Dutkay}
\address{Department of Mathematics\\
The University of Iowa\\
14 MacLean Hall\\
Iowa City, IA 52242-1419\\
U.S.A.}
\email{ddutkay@math.uiowa.edu}
\thanks{}
\subjclass{}
\keywords{}

\begin{abstract}
We give a complete description of spectrum of the wavelet Galerkin operator
$$R_{m_0,m_0}f(z)=\frac{1}{N}\sum_{w^N=z}\left|m_0\right|^2(w)f(w),\quad(z\in\mathbb{T})$$
associated to a a low-pass filter $m_0$ and a scale $N$, in the Banach spaces $C(\mathbb{T})$ and 
$L^{p}\left(\mathbb{T}\right)$, $1\leq p\leq\infty$.

\end{abstract}\maketitle
\tableofcontents
\section{\label{Intro}Introduction}
We begin with a short motivation of our study. For more background on wavelets and their connection to the wavelet 
Galerkin operator we refer the reader to \cite{Dau92}, \cite{BraJo} or \cite{HeWe}.
The wavelet analysis studies functions $\psi\in\ltwor$ with the property that 
$$\left\{2^{\frac{j}{2}}\psi\left(2^jx-k\right)\,|\, j,k\in\mathbb{Z}\right\}$$
is an orthonormal basis for $\ltwor$. Such functions are called wavelets. The scale (2 here) can be also 
any integer $N\geq2$. One way to construct wavelets is by multiresolutions. A multiresolution is 
a nest of subspaces $(V_j)_{j\in\mathbb{Z}}$ of $\ltwor$ with the following properties:
\begin{enumerate}
\item 
$V_j\subset V_{j+1}$, for all $j\in\mathbb{Z}$;
\item
$f\in V_j$ if and only if $f(Nx)\in V_{j+1}, (j\in\mathbb{Z})$;
\item
$$\bigcap_{j\in\mathbb{Z}}V_j=\{0\};$$
\item
$$\overline{\bigcup_{j\in\mathbb{Z}}V_j}=\ltwor;$$
\item 
There exists a function $\varphi\in V_0$ such that $\{\varphi(x-k)\,|\, k\in\mathbb{Z}\}$ is an orthonormal basis 
for $V_0$.
\end{enumerate}
\par
To build such a multiresolution one needs the function $\varphi$ called scaling function (or father function or refinable 
function). The scaling function satisfies a scaling equation:
$$\frac{1}{\sqrt{N}}\varphi\left(\frac{x}{N}\right)=\sum_{k\in\mathbb{Z}}a_k\varphi(x-k),\quad(x\in\mathbb{R}),$$
$a_k$ being some complex coefficients. The Fourier transform of the scaling equation is:
$$\sqrt{N}\widehat{\varphi}(N\xi)=m_0(\xi)\widehat{\varphi}(\xi),\quad(\xi\in\mathbb{R}),$$
where $m_0(\xi)=\sum_{k\in\mathbb{Z}}a_ke^{-ik\xi}$ is a $2\pi$-periodic function called low-pass filter. 
\par
Thus, the scaling functions $\varphi$ are determined by the low-pass filters $m_0$ and the construction of scaling functions 
has the low-pass filters as the starting point. 
\par
The multiresolution theory has shown that many of the properties of the scaling function $\varphi$ can be expressed in terms 
of the wavelet Galerkin operator associated to the filter $m_0$: 
$$R_{m_0,m_0}f(z)=\frac{1}{N}\sum_{w^N=z}\left|m_0\right|^2(w)f(w),\quad(z\in\mathbb{T}).$$
$\mathbb{T}$ is the unit circle, $f$ is some measurable function on $\mathbb{T}$, and we will identify 
functions on $\mathbb{T}$ with $2\pi$-periodic functions on $\mathbb{R}$. 
\par
For example, one needs the integer translates of the scaling function $\varphi(x-k)$, $k\in\mathbb{Z}$, to be orthonormal. To 
obtain this, a neccesary condition is the quadrature mirror filter condition:
$$\frac{1}{N}\sum_{w^N=z}\left|m_0\right|^2(w)=1,\quad(z\in\mathbb{T}),$$
which can be rewritten as $R_{m_0,m_0}1=1$. In \cite{Law91a} it is proved that the integer translates of the scaling function 
form an orthonormal set if and only if the constants are the only continuous functions that satisfy $R_{m_0,m_0}h=h$. So 
1 has to be a simple eigenvalue for the operator $R_{m_0,m_0}:C(\mathbb{T})\rightarrow C(\mathbb{T})$. Also, the regularity of the 
scaling function can be determined by the spectrum of $R_{m_0,m_0}$ (see \cite{Str96},\cite{RoSh}).
\par
We will impose some restrictions on $m_0$, restrictions that are custom in the setting of compactly supported wavelets:
\begin{equation}
m_0\mbox{ is a Lipschitz function};\\ \label{eq16}
\end{equation}
\begin{equation}
m_0\mbox{ has only a finite number of zeroes};\\ \label{eq17}
\end{equation}
\begin{equation}
m_0(0)=\sqrt{N};\\ \label{eq18} 
\end{equation}
\begin{equation}
R_{m_0,m_0}1=1. \label{eq19}
\end{equation}
In fact, for compactly supported wavelets, $m_0$ is a trigonometric polynomial, but for our purpose 
we can assume more generally that $m_0$ is Lipschitz. 
\par
The wavelet Galerkin operator $R_{m_0,m_0}$ bears several other names in the literature. It is also called the 
Ruelle operator because there are connections with the Ruelle-Perron-Frobenius theory for positive 
operators (see\cite{Bal00}), or transfer operator. We will use these names in the sequel. 
\par
 An extensive study of the spectral properties of the Ruelle operator can be found in \cite{BraJo}. We will gather 
some results from \cite{BraJo},\cite{Dutb} and add some new ones to give a complete picture of the spectrum of this 
Ruelle operator in the Banach spaces $C(\mathbb{T})$ and $L^p(\mathbb{T})$, $1\leq p\leq\infty$, answering in this way 
some questions posed in \cite{BraJo}.

\section{\label{Spectrum}The Spectrum of $R_{m_0,m_0}$}
In this section we present the results. We consider an integer $N\geq2$ and a function $m_0$ on $\mathbb{T}$ that satisfies (\ref{eq16})-(\ref{eq19}).
To $m_0$ we associate the Ruelle operator $R_{m_0,m_0}$ defined by 
$$R_{m_0,m_0}f(z)=\frac{1}{N}\sum_{w^N=z}\left|m_0\right|^2(w)f(w),\quad(z\in\mathbb{T}),$$
where $f$ is a measurable function on $\mathbb{T}$. We will see that $R_{m_0,m_0}$ is an operator on the spaces 
$C(\mathbb{T})$, $L^p(\mathbb{T})$ where $1\leq p\leq\infty$, and we will describe the spectrum and the eigenvalue spectrum 
of this operator on these spaces. 
\par
Before we give the results, some definitions and notations are needed. We denote by $R=R_{m_0,m_0}$. For a function 
$f$ on $\mathbb{T}$ and $\rho\in\mathbb{T}$ 
$$\alpha_{\rho}(f)(z)=f(\rho z),\quad(z\in\mathbb{T}).$$
For $\varphi\in\loner$, 
$$\operatorname*{Per}(\varphi)(x)=\sum_{k\in\mathbb{Z}}\varphi(x+2k\pi),\quad(x\in\mathbb{R}).$$
\par
We call a set $\{z_1,...,z_p\}$ a cycle of length $p$, and denote this by $z_1\rightarrow...\rightarrow z_p\rightarrow z_1$, 
if $z_1^N=z_2, z_2^N=z_3,...,z_{p-1}^N=z_p, z_p^N=z_1$ and the points $z_1,...,z_p$ are distinct. 
We call $z_1\rightarrow...\rightarrow z_p\rightarrow z_1$ an $m_0$-cycle if $|m_0|(z_i)=\sqrt{N}$ for $i\in\{1,...,p\}$.
\par
For a complex function $f$ on $\mathbb{T}$ and a positive integer $n$, 
$$f^{(n)}(z)=f(z)f\left(z^N\right)...f\left(z^{N^{n-1}}\right),\quad(z\in\mathbb{T}).$$

\begin{theorem}[The spectrum of $R$ on $C(T)$]\label{th21}
Let $m_0$ be a function satisfying (\ref{eq16})-(\ref{eq19}).
\begin{enumerate}
\item
The spectral radius of the operator $R:C(\mathbb{T})\rightarrow C(\mathbb{T})$ is equal to 1;
\item
Each point $\lambda\in\mathbb{C}$ with $|\lambda|<1$ is an eigenvalue for $R$, having infinite multiplicity and the spectrum 
of $R$ on $C(\mathbb{T})$ is the unit disk $\{\lambda\in\mathbb{C}\,|\, |\lambda|\leq1\}$;
\item
(The peripheral spectrum)
Let $C_1,...,C_n$ be the $m_0$-cycles, 
$$C_i=z_{1i}\rightarrow...\rightarrow z_{p_ii}\rightarrow z_{1i},\quad(i\in\{1,...,n\}).$$
Let $\lambda\in\mathbb{C}$, $|\lambda|=1$. Then $\lambda$ is an eigenvalue for $R$ if and only if $\lambda^{p_i}=1$ for some 
$i\in\{1,...,n\}$. The multiplicity of such a $\lambda$ equals the cardinality of the set 
$$\{i\in\{1,...,n\}\,|\,\lambda^{p_i}=1\}.$$
A basis for the eigenspace corresponding to $\lambda$ is obtained as follows:
\par
For $i\in\{1,...,n\}$ and $k\in\{1,...,p_i\}$ define
$$\varphi_{ki}(x)=\prod_{l=1}^{\infty}\frac{e^{-i\theta_i}\alpha_{z_{ki}}\left(m_0^{(p_i)}\right)\left(\frac{x}{N^{lp_i}}\right)}{\sqrt{N^{p_i}}},\quad(x\in\mathbb{R}),$$
where
$$e^{i\theta_i}=\frac{m_0(z_{1i})}{|m_0|(z_{1i})}...\frac{m_0(z_{p_ii})}{|m_0|(z_{p_ii})}.$$
Then 
$$g_{ki}=\alpha_{z_{ki}^{-1}}\left(\operatorname*{Per}|\varphi_{ik}|^2\right).$$
The basis for the eigenspace corresponding to the eigenvalue $\lambda$ is 
$$\{\sum_{k=1}^{p_i}\lambda^{-k+1}g_{ki}\,|\,i\in\{1,...,n\}\mbox{ with }\lambda^{p_i}=1\}.$$
Moreover, the functions in this basis are Lipschitz (or trigonometric polynomials when $m_0$ is one).
\end{enumerate}
\end{theorem}

\begin{proof}
(i) Take $f\in C(\mathbb{T})$ and $z\in\mathbb{T}$.
\begin{align*}
|Rf(z)|&\leq\frac{1}{N}\sum_{w^N=z}|m_0(w)|^2|f(w)|\\
&\leq\left\|f\right\|_{\infty}\frac{1}{N}\sum_{w^N=z}|m_0(w)|^2=\left\|f\right\|_{\infty}
\end{align*}
Therefore $\left\|Rf\right\|_{\infty}\leq\left\|f\right\|_{\infty}$ so the spectral radius is less then 1. But condition (\ref{eq19})
implies that 1 is an eigenvalue for $R$ so the spectral radius is 1.
\par
(ii) We begin with a lemma
\begin{lemma}\label{lema22}
If $z_0\rightarrow...\rightarrow z_{p-1}\rightarrow z_0$ is a cycle with $p$ large enough, then 
there exists a continuous function $f\neq0$ with $Rf=0$, such that $f(z_0)=1$, $f(z_i)=0$ for 
$i\in\{1,...,p-1\}$. 
\end{lemma}
\begin{proof}
To be able to produce such a function, we will need some conditions on the cycle. We will need $z_0$ and
$e^{\frac{2\pi i}{N}}z_0$ to be outside the set of zeroes of $m_0$. Because $m_0$ has only finitely many zeros,
this can be achieved as long as $p$ is big enough. We will also need $e^{\frac{2\pi i}{N}}z_0\neq z_l$ 
for $l\in\{1,...,p-1\}$, but this is true because, otherwise, $z_1=\left(e^{\frac{2\pi i}{N}}z_0\right)^N=z_l^N=z_{l+1}$ for some 
$l\in\{1,...,p-1\}$. 
\par 
So, when the cycle is long enough we have that $z_0,e^{\frac{2\pi i}{N}}z_0$ are outside the set of zeroes of $m_0$ 
and also $e^{\frac{2\pi i}{N}}z_0\neq z_l$ for all $l\in\{1,...,p-1\}$. Then we can choose a small interval $[a,b]$ 
(on $\mathbb{T}$) centered at $z_0$, such that 
\begin{equation}\label{eq21}
[a,b]\cup[e^{\frac{2\pi i}{N}}a,e^{\frac{2\pi i}{N}}b]\mbox{ contains no zeroes of }m_0;
\end{equation}
\begin{equation}\label{eq22}
[a,b]\cup[e^{\frac{2\pi i}{N}}a,e^{\frac{2\pi i}{N}}b]\mbox{ contains no }z_l, l\in\{1,...,p-1\};
\end{equation}
\begin{equation}\label{eq23}
\mbox{The intervals }[e^{\frac{2k\pi}{N}i}a,e^{\frac{2k\pi}{N}i}b], k\in\{0,...,N-1\}\mbox{ are disjoint}.
\end{equation}
\par
Define $f$ on $[a,b]$ continuously, to be 1 at $z_0$ and 0 at $a$ and $b$. Define $f$ on $[e^{\frac{2\pi i}{N}}a,e^{\frac{2\pi i}{N}}b]$ 
by 
$$f(z)=-\frac{1}{|m_0(z)|^2}|m_0|^2\left(e^{\frac{-2\pi i}{N}}z\right)f\left(e^{\frac{-2\pi i}{N}}z\right),\quad(z\in [e^{\frac{2\pi i}{N}}a,e^{\frac{2\pi i}{N}}b])$$
and define $f$ to be 0 everywhere else. $f$ is well defined because of (\ref{eq21}) and (\ref{eq23}). $f$ is continuous 
because it is $0$ at $a, b, e^{\frac{2\pi i}{N}}a$ and $e^{\frac{2\pi i}{N}}b$. It is also clear that $f(z_0)=1$ and 
$f(z_i)=0$ for $i\in\{1,...,p-1\}$ due to (\ref{eq22}).
\par
Now we check that $Rf=0$ which amounts to verifying that 
\begin{equation}\label{eq24}
\sum_{k=0}^{N-1}|m_0|^2\left(e^{\frac{2k\pi i}{N}}z\right)f\left(e^{\frac{2k\pi i}{N}}z\right)=0,\quad(z\in\mathbb{T})
\end{equation}
The only interesting case is when for some $k$, 
$$e^{\frac{2k\pi i}{N}}z\in[a,b]\cup[e^{\frac{2\pi i}{N}}a,e^{\frac{2\pi i}{N}}b].$$
So assume $e^{\frac{2k\pi i}{N}}z\in[a,b]$ for some $k\in\{0,...,N-1\}$. Then 
$$e^{\frac{2(k+1)\pi i}{N}}z\in[e^{\frac{2\pi i}{N}}a,e^{\frac{2\pi i}{N}}b]$$
and, using (\ref{eq23}), $f\left(e^{\frac{2l\pi i}{N}}z\right)=0$ for $l\in\{0,...,N-1\}\setminus\{k,k+1\}$.
(We use here notation modulo $N$ that is $N=0,N+1=1$ etc.)
Having theese, (\ref{eq24}) follows from the definition of $f$. 
\par
If 
$$e^{\frac{2k\pi i}{N}}z\in[e^{\frac{2\pi i}{N}}a,e^{\frac{2\pi i}{N}}b]$$
then 
$$e^{\frac{2(k-1)\pi i}{N}}\in[a,b]$$
and we can use the same argument as before to obtain (\ref{eq24}). This concludes the proof of the lemma.
\end{proof} 
We return to the prof of our theorem. Take $\lambda\in\mathbb{C}$ with $|\lambda|<1$. Choose a long enough cycle 
$z_0\rightarrow...\rightarrow z_{p-1}\rightarrow z_0$. Lemma (\ref{lema22}) produces a function $f_{z_0}\in C(\mathbb{T})$ 
with $Rf_{z_0}=0$, $f_{z_0}(z_i)=\delta_{0i}$ for $i\in\{0,...,p-1\}$. 
\par
Define
$$h_{z_0}(z)=\sum_{n=0}^{\infty}\lambda^n f\left(z^{N^n}\right),\quad(z\in\mathbb{T}).$$
(For $\lambda=0$ we make the convention $\lambda^0=1$.)
\par
The series is uniformly convergent because $\left\|f_{z_0}\left(z^{N^n}\right)\right\|_{\infty}=\left\|f\right\|_{\infty}$ 
for all $n\geq0$ and $|\lambda|<1$, so $h_{z_0}$ is continuous. 
\par
Also, if we use the fact that $R\left(f\left(z^{N^n}\right)\right)=f\left(z^{N^{n-1}}\right)$ for $n\geq1$, which is a consequence of the 
definition of $R$ and (\ref{eq19}), we have:
\begin{align*}
Rh_{z_0}(z)&=Rf_{z_0}(z)+\sum_{n=1}^{\infty}\lambda^nR\left(f\left(z^{N^n}\right)\right)\\
&=\lambda\sum_{n=1}^{\infty}\lambda^{n-1}f\left(z^{N^{n-1}}\right)=\lambda h_{z_0}
\end{align*}
\par
We evaluate $h_{z_0}$ at the points of the cycle $z_0,z_1,...,z_{p-1}$. Note that 
$$f_{z_0}\left(z_i^{N^n}\right)=f_{z_0}(z_{n+i})=
\left\{\begin{array}{rcl}
1&\mbox{ for }&n+i=0\mod p\\
0&\mbox{ otherwise }&
\end{array}\right.$$
(again, we use notation mod $p$, $z_{p}=z_0, z_{p+1}=z_1$, etc.)
\par
Hence,
$$h_{z_0}(z_0)=\sum_{m=0}^{\infty}\lambda^{mp}=\frac{1}{1-\lambda^p},$$
$$h_{z_0}(z_i)=\sum_{m=0}^{\infty}\lambda^{p-i+mp}=\frac{\lambda^{p-i}}{1-\lambda^p},\quad(i\in\{1,...,p-1\}),$$
so
$$(h_{z_0}(z_0),...,h_{z_0}(z_{p-1}))=\frac{1}{1-\lambda^p}(1,\lambda^{p-1},\lambda^{p-2},...,\lambda^2,\lambda).$$
\par
Now we make the same construction but considering the cycle starting from $z_k$. We obtain a function $f_{z_k}\in C(\mathbb{T})$ 
satisfying $Rf_{z_k}=0$, $f_{z_k}(z_i)=\delta_{ki}$ and
$$h_{z_k}(z)=\sum_{n=0}^{\infty}\lambda^nf_{z_k}\left(z^{N^n}\right)$$
has the properties $h_{z_k}\in C(\mathbb{T})$, $Rh_{z_k}=\lambda h_{z_k}$ and, for example, for $k=1$ we have the vector
$$(h_{z_1}(z_0),...,h_{z_1}(z_{p-1}))=\frac{1}{1-\lambda^p}(\lambda,1,\lambda^{p-1},\lambda^{p-2},...,\lambda^2).$$
Note that this vector is obtained from the previous one (the one corresponding to $z_0$), after a cyclic permutation. 
In fact the matrix 
$$(1-\lambda^p)\left( \begin{array}{cccc}
h_{z_0}(z_0) & h_{z_0}(z_1) & ... & h_{z_0}(z_{p-1}) \\
h_{z_1}(z_0) & h_{z_1}(z_1) & ... & h_{z_1}(z_{p-1}) \\
\vdots       & \vdots       &     & \vdots\\
h_{z_{p-1}}(z_0) & h_{z_{p-1}}(z_1) & ... & h_{z_{p-1}}(z_{p-1})
\end{array} \right)$$
is equal to 
$$\left( \begin{array}{ccccccc}
1&\lambda^{p-1}&\lambda^{p-2}&...&\lambda^{3}&\lambda^{2}&\lambda\\
\lambda&1&\lambda^{p-1}&...&\lambda^4&\lambda^3&\lambda^2\\
\vdots&\vdots&\vdots& &\vdots&\vdots&\vdots\\
\lambda^{p-2}&\lambda^{p-3}&\lambda^{p-4}&...&\lambda&1&\lambda^{p-1}\\
\lambda^{p-1}&\lambda^{p-2}&\lambda^{p-3}&...&\lambda^2&\lambda&1
\end{array}\right)
$$
Our goal is to prove that $h_{z_k}$ are linearly independent. We can achieve this if we show that the matrix is nonsingular. For this,
look at the entries below the diagonal. We note that, on each row, the part below the diagonal can be obtained from the previous row times 
$\lambda$. Therefore, if we subtract from the $p-1$-st row $\lambda$ times the $p-2$-nd row, substract from the $p-2$-nd row 
$\lambda$ times the $p-3$-rd row,..., substract from the $1$-st row $\lambda$ times the $0$-th row, we obtain an upper triangular  
matrix having $1-\lambda^p$ on each diagonal entry and which has the same determinant as the initial one. Since $|\lambda|<1$, 
this matrix will be nonsingular so $h_{z_0},h_{z_1},...,h_{z_{p-1}}$ are linearly independent eigenvectors that correspond to the 
eigenvalue $\lambda$. As $p$ can be chosen as big as we want, the multiplicity of $\lambda$ is infinite.
\par
(iii) See \cite{Dutb}.
\end{proof}

\begin{theorem}[The spectrum of $R$ on $\linfr$]\label{th23}
Let $m_0$ be a function on $\mathbb{T}$ satisfying (\ref{eq16})-(\ref{eq19}). 
\begin{enumerate}
\item
The spectral radius of the operator $R:\linfr\rightarrow\linfr$ is equal to 1 and the spectrum of $R$ 
is the unit disk $\{\lambda\in\mathbb{C}\,|\, |\lambda|\leq 1\}$.
\item
Each point $\lambda\in\mathbb{C}$ with $|\lambda|\leq1$ is an eigenvalue for $R$ of infinite multiplicity.
\end{enumerate}
\end{theorem}

\begin{proof}
(i) The argument used in the proof of theorem \ref{th21} applies here to obtain the spectral radius equal to 1 and the fact 
that the spectrum is the unit disk will follow from (ii).
\par
(ii) If $|\lambda|<1$ then the assertion follows trivialy from theorem \ref{th21} (ii). It remains to consider the case 
$|\lambda|=1$. Define 
$$\varphi(x)=\prod_{n=1}^{\infty}\frac{m_0\left(\frac{x}{N^n}\right)}{\sqrt{N}},\quad(x\in\mathbb{R}).$$ 
$\varphi$ is a well defined, continuous function in $\ltwor$ and $\operatorname*{Per}|\varphi|^2$ is a Lipschitz 
function on $\mathbb{T}$ (see \cite{BraJo}). Also, $\operatorname*{Per}|\varphi|^2(0)=1$, $\varphi(0)=1$ and 
$$\varphi(x)=\frac{m_0\left(\frac{x}{N}\right)}{\sqrt{N}}\varphi\left(\frac{x}{N}\right),\quad(x\in\mathbb{R}),$$
($\varphi$ is the Fourier transform of a scaling function).
\par
Now, consider a function $f\in\linfr$ with the property that $f(x)=\frac{1}{\lambda}f\left(\frac{x}{N}\right)$ a.e. on $\mathbb{R}$, 
and take $h_f=\operatorname*{Per}\left(f|\varphi|^2\right)$. Clearly, $|h_f(z)|\leq\left\|f\right\|_{\infty}\operatorname*{Per}|\varphi|^2(z)$ for 
$z\in\mathbb{T}$ so $h_f$ is an $\linft$ function. We want to prove that 
$$Rh_f=\lambda h_f.$$
\par
We have
\begin{align*}
h_f(x)&=\sum_{k\in\mathbb{Z}}f(x+2k\pi)|\varphi|^2(x+2k\pi)\\
&=\sum_{k\in\mathbb{Z}}\frac{1}{\lambda}f\left(\frac{x+2k\pi}{N}\right)\frac{1}{N}|m_0|^2\left(\frac{x+2k\pi}{N}\right)|\varphi|^2\left(\frac{x+2k\pi}{N}\right)\\
&=\frac{1}{N}\frac{1}{\lambda}\sum_{l=0}^{N-1}\sum_{m\in\mathbb{Z}}f\left(\frac{x+2(Nm+l)\pi}{N}\right)|m_0|^2\left(\frac{x+2(Nm+l)\pi}{N}\right)\cdot\\
&\cdot|\varphi|^2\left(\frac{x+2(Nm+l)\pi}{N}\right)\\
&=\frac{1}{N}\frac{1}{\lambda}\sum_{l=0}^{N-1}|m_0|^2\left(\frac{x+2l\pi}{N}\right)\sum_{m\in\mathbb{Z}}f|\varphi|^2\left(\frac{x+2k\pi}{N}+2m\pi\right)\\
&=\frac{1}{\lambda}Rh_f.
\end{align*}
so, we have indeed $Rh_f=\lambda h_f$.
\par
Next, we argue why the vector space
$$\left\{h_f\,|\, f\in\linfr, f(x)=\frac{1}{\lambda}f\left(\frac{x}{N}\right)\mbox{ a.e. on }\mathbb{R}\right\}$$ 
is infinite dimensional.
\par
For this, we prove first that if $h_f=0$ then $f=0$. Indeed, if $h_f=0$ then 
\begin{equation}\label{eq25}
f(x)|\varphi|^2(x)=-\sum_{k\in\mathbb{Z}\setminus\{0\}}f(x+2k\pi)|\varphi|^2(x+2k\pi).
\end{equation}
We claim that the term on the right converges to 0 as $x\rightarrow 0$. We have
$$
\left|\sum_{k\in\mathbb{Z}\setminus\{0\}}f(x+2k\pi)|\varphi|^2(x+2k\pi)\right|\leq$$
$$\leq\left\|f\right\|_{\infty}\sum_{k\in\mathbb{Z}\setminus\{0\}}|\varphi|^2(x+2k\pi)=\left\|f\right\|_{\infty}(\operatorname*{Per}|\varphi|^2(x)-|\varphi|^2(x))
\rightarrow 0,$$
because both $\operatorname*{Per}|\varphi|^2$ and $|\varphi|^2$ are continuous and their value at 0 is 1. Then, using (\ref{eq25}),
we obtain $f(x)\rightarrow 0$ as $x\rightarrow 0$. But we know that $f(x)=\frac{1}{\lambda}f\left(\frac{x}{N}\right)$ a.e. 
on $\mathbb{R}$. So for a.e. $x$ we have 
$$f\left(\frac{x}{N^n}\right)=\lambda^nf(x),\quad\mbox{ for all } n.$$
This implies that 
$$\left|f\left(\frac{x}{N^n}\right)\right|=|f(x)|,\quad(n\in\mathbb{N})$$
and, coupled with the limit of $f$ at $0$, it entails that $f$ is constant 0 almost everywhere. 
\par
Having these, we try to construct a set of $p$ linearly independent functions $h_f$ with $p\in\mathbb{N}$ arbitrary. 
Take $p$ linearly independent functions $g_1,...,g_p$ in $L^{\infty}([-N,-1]\cup[1,N])$. Define $f_i$, $i\in\{1,...,p\}$
 on $\mathbb{R}$ as follows: let $f_i(x)=g_i(x)$ on $[-N,-1]\cup[1,N]$ and extend it on $\mathbb{R}$ requiring that 
$$\frac{1}{\lambda}f_i\left(\frac{x}{N}\right)=f_i(x),\quad(x\in\mathbb{R}).$$
That is, for $x\in[-\frac{1}{N^l},-\frac{1}{N^{l+1}}]\cup[\frac{1}{N^{l+1}}\cup\frac{1}{N^l}]$ 
$$f_i(x)=\lambda^{l+1}g_i(N^{l+1}x),$$
for all $l\in\mathbb{Z}$. Since $|\lambda|=1$, $f_1,...,f_p$ are in $\linfr$ and they are linearly independent because 
$g_1,...,g_p$ are. $h_{f_1},...,h_{f_p}$ are linearly independent by the following argument: if for some complex constants 
$a_1,...,a_p$ we have $a_1h_{f_1}+...a_ph_{f_p}=0$ then $h_{a_1f_1+...+a_pf_p}=0$ so $a_1f_1+...+a_pf_p=0$ and $a_1=...=a_p=0$ by linear 
independence. Since we proved that $Rh_{f_i}=\lambda_{f_i}$, $i\in\{1,...,p\}$, and since $p$ is arbitrary, it follows that 
the multiplicity of the eigenvalue $\lambda$ is infinite.
\end{proof}
\begin{theorem}[The spectrum of $R$ on $L^p(\mathbb{T})$]\label{th24}
Let $m_0$ be a function on $\mathbb{T}$ satisfying (\ref{eq16})-(\ref{eq19}) and $1\leq p<\infty$.
\begin{enumerate}
\item
The spectral radius of the operator $R:L^p(\mathbb{T})\rightarrow L^p(\mathbb{T})$ is equal to $N^{\frac{1}{p}}$ and the spectrum 
of $R$ is the disk $\{\lambda\in\mathbb{C}\,|\, |\lambda|\leq N^{\frac{1}{p}}\}$.
\item
Each point $\lambda\in\mathbb{C}$ with $|\lambda|<N^{\frac{1}{p}}$ is an eigenvalue for $R$ of infinite multiplicity.
\item
There are no eigenvalues of $R$ with $|\lambda|=N^{\frac{1}{p}}$.
\end{enumerate}
\end{theorem}
\begin{proof}
(i) is proved in \cite{BraJo} but we present here a different argument that we will need for (iii) also. Take $f\in L^p(\mathbb{T})$.
\begin{align*}
\left\|Rf\right\|_p&=\left(\int_0^{2\pi}\left|\frac{1}{N}\sum_{k=0}^{N-1}|m_0|^2f\left(\frac{\theta+2k\pi}{N}\right)\right|^p\,d\theta\right)^{\frac{1}{p}}\\
&\leq\left(\int_0^{2\pi}\left(\frac{1}{N}\sum_{k=0}^{N-1}|m_0|^2|f|\left(\frac{\theta+2k\pi}{N}\right)\right)^p\,d\theta\right)^{\frac{1}{p}}
\end{align*}
Since 
$$\frac{1}{N}\sum_{k=0}^{N-1}|m_0|^2\left(\frac{\theta+2k\pi}{N}\right)=1,\quad(\theta\in[0,2\pi])$$
and $x\mapsto x^p$ is convex, we can use Jensen's inequality: 
\begin{align*}
\left(\frac{1}{N}\sum_{k=0}^{N-1}|m_0|^2|f|\left(\frac{\theta+2k\pi}{N}\right)\right)^p&\leq\frac{1}{N}\sum_{k=0}^{N-1}|m_0|^2|f|^p\left(\frac{\theta+2k\pi}{N}\right)\\
&\leq\sum_{k=0}^{N-1}|f|^p\left(\frac{\theta+2k\pi}{N}\right)
\end{align*}
For the last inequality we used the fact that $|m_0|^2\leq N$ which follows from (\ref{eq19}).
\par
Also, by a change of variable, 
\begin{align*}
\left(\int_0^{2\pi}\sum_{k=0}^{N-1}|f|^p\left(\frac{\theta+2k\pi}{N}\right)\,d\theta\right)^{\frac{1}{p}}&=
\left(\sum_{k=0}^{N-1}N\int_{\frac{2k\pi}{N}}^{\frac{2(k+1)\pi}{N}}|f(\theta)|^p\,d\theta\right)^\frac{1}{p}\\
&=N^{\frac{1}{p}}\int_0^{2\pi}|f(\theta)|^p\,d\theta.
\end{align*}
Putting together the previous equalities and inequalities we obtain that $\left\|Rf\right\|_p\leq N^{\frac{1}{p}}\left\|f\right\|_p$. 
This implies that the norm and the spectral radius of the operator $R:L^p(\mathbb{T})\rightarrow L^p(\mathbb{T})$ are less than $N^{\frac{1}{p}}$. 
A result of R. Nussbaum (see \cite{BraJo}) shows that every $\lambda\in\mathbb{C}$ with $1<|\lambda|<N^{\frac{1}{p}}$ is 
an eigenvalue of $R$ of infinite multiplicity. Also theorem \ref{th23} shows that all $\lambda\in\mathbb{C}$ with $|\lambda|\leq 1$ 
is an eigenvalue of $R$ of infinite multiplicity. This establishes (i) and (ii). 
\par
It remains to prove that (iii) is valid. Suppose there is a function $f\in L^p(\mathbb{T})$ and $\lambda\in\mathbb{C}$ 
such that $|\lambda|=N^{\frac{1}{p}}$ and $Rf=\lambda f$.
Then $\left\|Rf\right\|_p=N^{\frac{1}{p}}\left\|f\right\|_p$ so 
we have equalities in all inequalities that we used for proving (i). In particular, we have
$$\int_0^{2\pi}\frac{1}{N}\sum_{k=0}^{N-1}|m_0|^2|f|^p\left(\frac{\theta+2k\pi}{N}\right)\,d\theta=\int_0^{2\pi}\sum_{k=0}^{N-1}|f|^p\left(\frac{\theta+2k\pi}{N}\right)\,d\theta$$
and, since $\frac{|m_0|^2}{N}\leq1$ the corresponding terms of the sums must be equal: for $k\in\{0,...,N-1\}$, 
$$\int_0^{2\pi}\frac{|m_0|^2\left(\frac{\theta+2k\pi}{N}\right)}{N}|f|^p\left(\frac{\theta+2k\pi}{N}\right)\,d\theta=
\int_0^{2\pi}|f|^p\left(\frac{\theta+2k\pi}{N}\right)\,d\theta$$
Therefore, utilizing again $|m_0|^2\leq N$,
$$\frac{|m_0|^2\left(\frac{\theta+2k\pi}{N}\right)}{N}|f|^p\left(\frac{\theta+2k\pi}{N}\right)=
|f|^p\left(\frac{\theta+2k\pi}{N}\right)$$
for almost every $\theta\in[0,2\pi]$ and for all $k\in\{0,...,N-1\}$. But this implies that $\frac{1}{N}|m_0|^2|f|^p=|f|^p$ almost everywhere on 
$\mathbb{T}$. However, $m_0$ is continuous and has finitely many zeroes and, because $\sum_{w^N=z}|m_0|^2(w)=N$ for all 
$z\in\mathbb{T}$, this implies that $|m_0|^2(z)=N$ for at most finitely many points so $f$ must be 0 almost everywhere. In conclusion,
there are no eigenvalues $\lambda$ of modulus $N^{\frac{1}{p}}$ and the proof of the theorem is complete.
\end{proof}

\end{document}